\documentclass[conference, nofonttune,9pt]{IEEEtran}
\usepackage[normalem]{ulem}
\usepackage{amsmath}
\usepackage{amsfonts}
\usepackage{hyperref}
\usepackage{graphicx}
\usepackage{cleveref}
\usepackage{algorithm, algpseudocode,subcaption}
\usepackage{xcolor}
\usepackage{tikz}
\usepackage{caption}
\usepackage{listings}
\usepackage{array}
\usepackage{booktabs}
\usepackage{diagbox}
\usepackage{float}
\usepackage{geometry}
\newtheorem{theorem}{Theorem}

\newtheorem{definition}{Definition}

\title{\title{Scalable Binary CUR Low-Rank Approximation Algorithm}}
\author{Bowen Su }
\definecolor{officegreen}{rgb}{0.0, 0.5, 0.0}
\begin{document}

\maketitle

\section{Introduction}
Large-scale matrices from applications such as social networks~\cite{liben2003link} and genomics~\cite{alter2000singular} present significant challenges in modern data analysis. These matrices often exhibit low-rank structures, which can be exploited through matrix factorization techniques to derive compact and informative representations. Classical methods for low-rank matrix approximation, including Singular Value Decomposition (SVD)~\cite{golub1965calculating}, rank-revealing QR~\cite{gu1996efficient}, and rank-revealing LU decomposition~\cite{pan2000existence}, are well-established and provide reliable solutions. However, their computational costs scale rapidly with matrix size. Specifically, SVD has a complexity of \(\mathcal{O}(n^2 m)\) for a matrix \(A \in \mathbb{R}^{n \times m}\) with \(n \geq m\), while rank-revealing QR and LU decomposition require \(\mathcal{O}(n m^2)\). Such high computational demands render these methods impractical for applications involving extremely large datasets, creating significant challenges in terms of scalability and efficiency for large-scale data analysis. To address these limitations, CUR Low-Rank Approximaation~\cite{drineas2008cur,boutsidis2014optimal} has emerged as a potential alternative, offering efficient approximation of matrices with inherent low-rank structures.
\subsection{CUR Low-Rank Approximaation}
\begin{figure}[H]
    \centering
    \includegraphics[width=\linewidth]{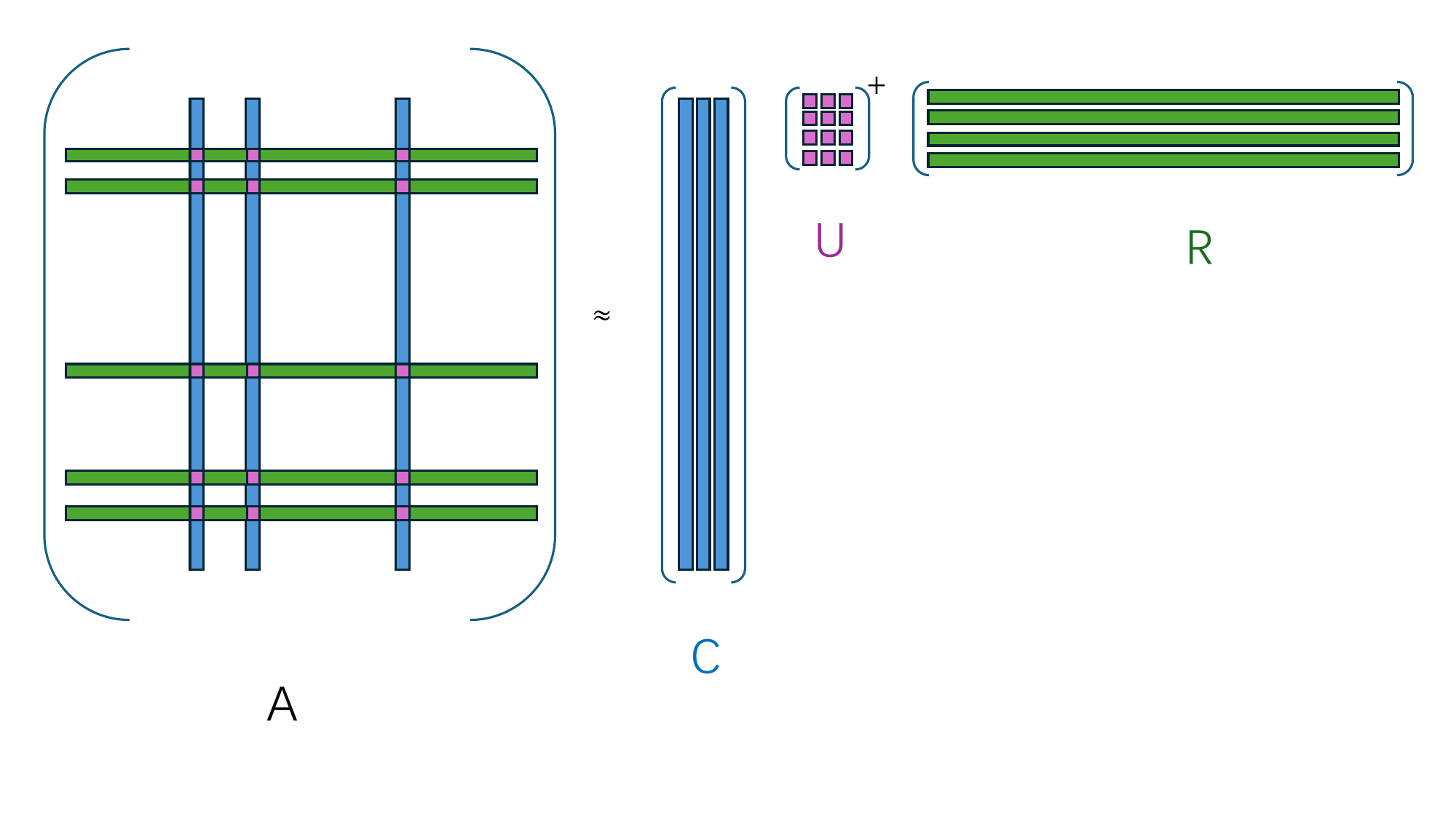}
    \caption{CUR Approximation}
    \label{fig:cur}
\end{figure}

Let \( A \in \mathbb{R}^{m \times n} \) be a given matrix, and consider two sets of indices \( I \subset [m] \) and \( J \subset [n] \), corresponding to selected rows and columns, respectively. Define \( A_{I,:} \) as the subset of rows in \( A \) indexed by \( I \), and \( A_{:,J} \) as the subset of columns indexed by \( J \). Suppose the submatrix \( A_{I,J} \), formed by the intersection of rows \( I \) and columns \( J \), is nonsingular. In this case, the approximation  
\[
A \approx CUR,
\]  
where \(C = A(:,J), U= A^{\dagger}(I,J), R = A(I,:), \)
is referred to as a CUR approximiation of \( A \). For a visual representation, see \Cref{fig:cur}.
There has been extensive research on CUR Low-Rank Approximations for matrices. For readers interested in further details and additional references, please see the following works \cite{aldroubi2019cur,allen2024maximal,anderson2015spectral,anderson2017efficient,boutsidis2014optimal,CaiCUR2021_2,caiCUR2021,cao2024randomized,chaturantabut2010nonlinear,chen2020efficient,civril2009selecting,dong2021simpler,drineas2008relative,drmac2016new,fornace2024column,georgieva2017best,gidisu2024restricted,goreinov2010find,goreinov2001maximal,goreinov2011quasioptimality,goreinov1997theory,gu1996efficient,hamm2019cur,HAMMCUR2020,KeatonCUR,HammCUR2021_2,ida2020fast,mahoney2009cur,park2024accuracy,park2024low,sorensen2016deim,thurau2012deterministic,voronin2017efficient,wang2012scalable,wang2013improving,wang2016towards,xia2024,zamarashkin2016new} and the references cited therein. There are different strategies to obtain such CUR approximations or decomposition.  {This paper explore maximal volume optimization strategy within the framework of parallel computing.} 
\subsection{Maximal Volume Optimization}
\begin{definition}
Given a square matrix $A\in \mathbb{R}^{r\times r},$ the volume of it is defined as absolute value of its determinant, i.e.,
\[
\operatorname{vol}(A) = |det(A)|.
\]
\end{definition}
The CUR approximation with \( r \) selected columns and rows is usually formulated as the following optimization problem \cite{allen2024maximal,civril2009selecting, goreinov2010find,goreinov2001maximal, goreinov2011quasioptimality,mikhalev2018rectangular,schneider2010error,zamarashkin2016new}  
\[
\operatorname*{argmax}_{\substack{|\hat{I}| = |\hat{J}| = r}} \operatorname{vol}(A(\hat{I}, \hat{J})).
\]
The aforementioned optimization problem has numerous practical applications such as identifying the optimal nodes for polynomial interpolation over a given domain on appropriately discretized meshes~\cite{sommariva2009computing} and preconditioning the iterative solution of linear least-squares problems~\cite{arioli2015preconditioning}.

\begin{theorem} [cf. \cite{goreinov2001maximal}]\label{thm:error_estimate}
Let \( A \) be a matrix of size \( m \times n \) and \( r\) be an integer satisfying \( 1 \leq r < \min\{m, n\}) \). Suppose that \( I \subset \{1, \cdots, m\} \) with \( \#(I) = r \) and \( J \subset \{1, \cdots, n\} \) with \( \#(J) = r \) such that \( A_{I,J} \) has the maximal volume among all \( r \times r \) submatrices of \( A \). Then, the Chebyshev norm of the residual matrix satisfies the following estimate:
\begin{equation}
\label{skapp}
\|A - A_r\|_C \le (r + 1)\sigma_{r+1}(A),
\end{equation}
where $A_r= A_{:,J} A_{I,J}^{\dagger} A_{I,:}$ is the $r$-cross approximation of $A$, 
$\sigma_{r+1}(A)$ is the $(r+1)^{th}$ singular value of $A$.
\end{theorem}
Recently, under the assumption in \Cref{thm:error_estimate}, Allen et.al \cite{allen2024maximal} improve the error estimate bounds to the following the results.
\begin{theorem}[cf. \cite{allen2024maximal}]\label{thm:improved_error_estimate}
Under the assumption in \Cref{thm:error_estimate}, there is following estimate: 
\[
\|A - A_{r}\|_{C} \leq \frac{(r+1)\sigma_{r+1}(A)}{\sqrt{1+\sum\limits_{k=1}^r \frac{\sigma_{r+1}^2(A)}{\sigma_{k}^{2}(A)}}}
\]
\end{theorem}

\Cref{thm:error_estimate} and \Cref{thm:improved_error_estimate} tells us that the volume maximization can lead to a quasi-best low-rank approximation. However, identifying maximum-volume submatrices presents considerable computational challenges, as the problem is proven to be NP-hard~\cite{bartholdi1982good}. This complexity has driven the development and adoption of greedy method-adaptive cross approximation-to provide practical solutions~\cite{bebendorf2000approximation}. At its core, this method is highly connected to incomplete LU decomposition of the matrix~\cite{cortinovis2020maximum}. We firstly introduce CUR Low-Rank Approximaation via Adaptive Cross Approximation Algorithm (ACA)~\cite{bebendorf2000approximation} as follows.
\begin{algorithm}[H]
\caption{Scalable Binary CUR Low-Rank Approximation Algorithm}
\label{alg:aca}
\begin{algorithmic}[1]
\State \textbf{Input:} Matrix $A \in \mathbb{R}^{m\times n}$, target rank $r$.
\State{Initialize \( R_0 := A \), \( I := \{\}, J := \{\}. \)}
\For{\( k = 0, \ldots, m-1 \)}

\State \(
    (i_{k+1}, j_{k+1}) := \operatorname*{argmax}_{i,j} |R_k(i,j)|.
    \)
\State \( I \gets I \cup \{i_{k+1}\}, J \gets J \cup \{j_{k+1}\}. \)
\State\(
R_{k+1} := R_k - \frac{1}{R_k(i_{k+1}, j_{k+1})} R_k(:, j_{k+1}) R_k(i_{k+1}, :).
    \)
\EndFor
\State \textbf{Output: \(C = A(:,J), U = A(I,J), R = A(I,:). \)}
\end{algorithmic}
\end{algorithm}
\noindent
Each iteration in \Cref{alg:aca} identifies a single pivot by searching for the largest absolute value in the residual matrix. Although local maxima searches can be conducted in parallel, a global reduction is necessary to determine the unique pivot element, thus creating an unavoidable sequential bottleneck. Once the pivot is found, the algorithm performs a rank-1  update on the residual, and this updated matrix in turn dictates the location and value of the next pivot. Consequently, there is a strict iteration-by-iteration data dependency, preventing the selection of multiple pivots simultaneously. 
In this paper, we propose the Scalable Binary CUR Low-Rank Approximation Algorithm. This framework is specifically designed to explore scalable CUR approximations for large-scale matrices by selecting rows and columns in parallel.

\section{Scalable Binary CUR Low-Rank Approximaation Algorithm}
\subsection{Overview}
The proposed \emph{Scalable Binary CUR Low-Rank Approximation Algorithm} efficiently approximates a large-scale matrix \( A \in \mathbb{R}^{n \times m} \) by decomposing it into three smaller matrices, \( C \), \( R \), and \( U \). First, a binary parallel selection process, as detailed in \Cref{alg:rowSelection}, identifies representative subsets of rows and columns from \( A \). The matrix \( C \) is constructed by extracting the columns of \( A \) corresponding to the selected indices, while \( R \) is formed by extracting the rows based on these indices. Finally, \( U \) is computed as the pseudoinverse of the intersection submatrix \( A(\text{rows}, \text{cols}) \). Together, these three matrices yield an efficient low-rank approximation of \( A \).

\begin{algorithm}[H]
\caption{Scalable Binary CUR Low-Rank Approximaation Algorithm}
\label{alg:Algorithm_CUR}
\begin{algorithmic}[1]
\State \textbf{Input:} Matrix \(A \in \mathbb{R}^{n \times m}\), selected row count \(r\), selected column count \(c\)
\State \textcolor{officegreen}{// Perform the following steps in parallel}
\State Use \Cref{alg:rowSelection} to determine the set of row indices \(I\).
\State Use \Cref{alg:rowSelection} to determine the set of column indices \(J\).

\State \textbf{Output: \(C = A(:,J), U = A(I,J), R = A(I,:). \)}
\end{algorithmic}
\end{algorithm}
\begin{algorithm}[!]
\caption{Blockwise Adaptive Cross Approximation}
\label{alg:rowSelection}
\begin{algorithmic}[1]
\State \textbf{Input:} Matrix \(A \in \mathbb{R}^{n \times m}\), number of selected indices \(r\), number of blocks \(b\), \texttt{axis} (0 for rows, 1 for columns)
\State \textbf{Output:} Selected indices $I$

\State $I = \{\}$  \Comment{Initialize selected indices }
\If{\texttt{axis} = 0} 
    \State Divide \(A\) into \(b\) blocks along rows: \(A_{\text{block}}^k \in \mathbb{R}^{n_k \times m}\)
    \State \texttt{dimBlockSize} $\gets n_k$, \texttt{dimSize} $\gets n$
\Else
    \State Divide \(A\) into \(b\) blocks along columns: \(A_{\text{block}}^k \in \mathbb{R}^{n \times m_k}\)
    \State \texttt{dimBlockSize} $\gets m_k$, \texttt{dimSize} $\gets m$
\EndIf

\For{$i = 0$ to $r - 1$} \Comment{Select \(r\) indices}
    \State \texttt{localMaxNorm}, \texttt{localMaxIdx} $\gets -1, -1$ 
    \State{\textcolor{officegreen}{// Blockwise parallel selection of rows/columns with maximal norm}}
\For{$i = 0$ to $b \cdot \texttt{dimBlockSize} - 1$ (parallel)} 
    \State \(k \gets \lfloor i / \texttt{dimBlockSize} \rfloor\) 
    \State \(j \gets i \bmod \texttt{dimBlockSize}\) 
    \State \( \texttt{globalIdx} \gets i \)

    \If{\texttt{axis} = 0}
        \State \( \texttt{localNorm} \gets \operatorname{norm}(A_{\text{block}}^{k}(j,:))\)
    \Else
        \State \( \texttt{localNorm} \gets \operatorname{norm}(A_{\text{block}}^{k}(:,j)) \)
    \EndIf

    \If{\texttt{localNorm} $>$ \texttt{localMaxNorm}}
        \State \texttt{localMaxNorm} $\gets$ \texttt{localNorm}
        \State \texttt{localMaxIdx} $\gets$ \texttt{globalIdx}
    \EndIf
\EndFor

   \For{$k = 0$ to $b - 1$ (parallel)} 
    \If{\texttt{localMaxNorm} $>$ \texttt{maxNorm}}
        \State \texttt{maxNorm} $\gets$ \texttt{localMaxNorm}
        \State \texttt{maxIdx} $\gets$ \texttt{localMaxIdx}
    \EndIf
    \EndFor
    \State $k_{\text{block}} \gets \texttt{maxIdx} \div \texttt{dimBlockSize}$
    \State $j_{\text{local}} \gets \texttt{maxIdx} \bmod \texttt{dimBlockSize}$
    \If{\texttt{axis} = 0}
        \State \texttt{sharedMaxVec} $\gets A_{\text{block}}^{k}(j_{\text{local}}, :)$
    \Else
        \State \texttt{sharedMaxVec} $\gets A_{\text{block}}^{k}(:, j_{\text{local}})$ 
    \EndIf
    \State{\textcolor{officegreen}{// Parallel Blockwise update using shared max vector}}     \For{$k = 0$ to $b - 1$ (parallel)}
        \If{\texttt{axis} = 0}
            \State \resizebox{0.7\linewidth}{!}{\(A_{\text{block}}^k \gets A_{\text{block}}^k - 
                    \frac{A_{\text{block}}^k \cdot\texttt{sharedMaxVec}^\top}{\texttt{maxNorm}} \cdot \texttt{sharedMaxVec}\)}
        \Else
            \State \resizebox{0.7\linewidth}{!}{\(A_{\text{block}}^k \gets A_{\text{block}}^k - 
                    \texttt{sharedMaxVec} \cdot 
                    \left(\frac{\texttt{sharedMaxVec}^\top \cdot A_{\text{block}}^k}{\texttt{maxNorm}}\right)\)}
        \EndIf
    \EndFor
    \State $I(i) \gets \texttt{maxIdx}$ 
\EndFor

\State \Return $I$
\end{algorithmic}
\end{algorithm}

\subsection{Blockwise Adaptive Cross Algorithm}
The Blockwise Adaptive Cross Approximation algorithm is designed to efficiently select a subset of \(r\) representative indices from a large matrix \(A \in \mathbb{R}^{n \times m}\) by iteratively extracting the most significant row or column vectors. The algorithm begins by partitioning the matrix into \(b\) blocks along a specified axis, determined by the parameter \texttt{axis} (where 0 indicates a row-wise partition and 1 indicates a column-wise partition). When \texttt{axis} is set to 0, \(A\) is divided into blocks \(A_{\text{block}}^k \in \mathbb{R}^{n_k \times m}\) along the rows; if \texttt{axis} is set to 1, the matrix is partitioned into blocks \(A_{\text{block}}^k \in \mathbb{R}^{n \times m_k}\) along the columns. In this context, the variable \texttt{dimBlockSize} refers to the number of rows or columns in each block, while \texttt{dimSize} represents the total number of rows \(n\) or columns \(m\) in the matrix.

An empty global index set \(I\) is then initialized to store the indices of the selected rows or columns. The algorithm proceeds over \(r\) iterations, each intended to identify one key index. In each iteration, a parallel search is conducted within each block to compute the Euclidean norm of every candidate vector (row or column, depending on the specified axis). Each block determines its local maximum norm and the corresponding local index. These local maxima are subsequently combined through a reduction step that identifies the global maximum norm and its associated index, denoted as \(\texttt{maxIdx}\).

Once the vector corresponding to \(\texttt{maxIdx}\) is identified, it is extracted from its block and referred to as \(\texttt{sharedMaxVec}\). This vector is then employed to perform a rank-1 update on every block of the matrix. Specifically, if the operation is row-wise (i.e., \texttt{axis} = 0), the update involves subtracting from each block the product of the block and the transpose of \(\texttt{sharedMaxVec}\) divided by the squared norm of \(\texttt{sharedMaxVec}\), multiplied by \(\texttt{sharedMaxVec}\) itself. Conversely, if the operation is column-wise (i.e., \texttt{axis} = 1), the subtraction is carried out by multiplying \(\texttt{sharedMaxVec}\) with the appropriate normalized factor computed from the block. This subtraction effectively removes the dominant component captured by the selected vector, ensuring that subsequent iterations focus on extracting new, non-redundant information.

 The complete pseudo code is provided in \Cref{alg:rowSelection} for detailed reference.

The complexity of the  Block Adaptive Cross Algorithm is determined by the operations performed in each iteration. Computing row norms for \(n_k\) rows within a block requires \(\mathcal{O}(n_k m)\) operations, and with \(b\) blocks, this totals \(\mathcal{O}\left(\frac{nm}{b}\right)\) per iteration. Updating the matrix, which includes projection, also requires \(\mathcal{O}\left(\frac{nm}{b}\right)\) operations per iteration. 

For \(r\) iterations, the total complexity of the algorithm is:
\(
\mathcal{O}\left(r\cdot\frac{nm}{b}\right).
\)
 This represents a significant improvement over the sequential Adaptive Cross Approximation algorithm’s complexity of \(O(r\cdot nm)\), particularly for large matrices and high levels of parallelism.
\section{Numerical Experiments}
All experiments are conducted in a computing environment powered by an AMD EPYC 7H12 64-Core Processor running at 2.6 GHz. For matrix operations, we utilized Armiddol\cite{sanderson2016armadillo,sanderson2019practical}, a high-performance matrix library designed to efficiently handle large-scale linear algebra computations. 

\subsection{Approximation Error Study of \Cref{alg:Algorithm_CUR}}
In this section, we evaluate the performance of \Cref{alg:Algorithm_CUR} by testing its approximation capabilities on two types of data: \textit{Hilbert matrices} and \textit{synthetic low-rank matrices}. The key metric used for evaluation is the Frobenius norm relative error, defined as:
\[
\text{Relative Error} = \frac{\|A - C U R\|_F}{\|A\|_F},
\]
where \(C\), \(U\), and \(R\) are the matrices resulting from the CUR Low-Rank Approximaation of \(A\). This metric provides a quantitative assessment of the quality of the approximation and is consistently applied across all experiments.

\paragraph{Experiments on Hilbert Matrices}

\begin{figure}[H]
\centering
\includegraphics[width=0.5\textwidth]{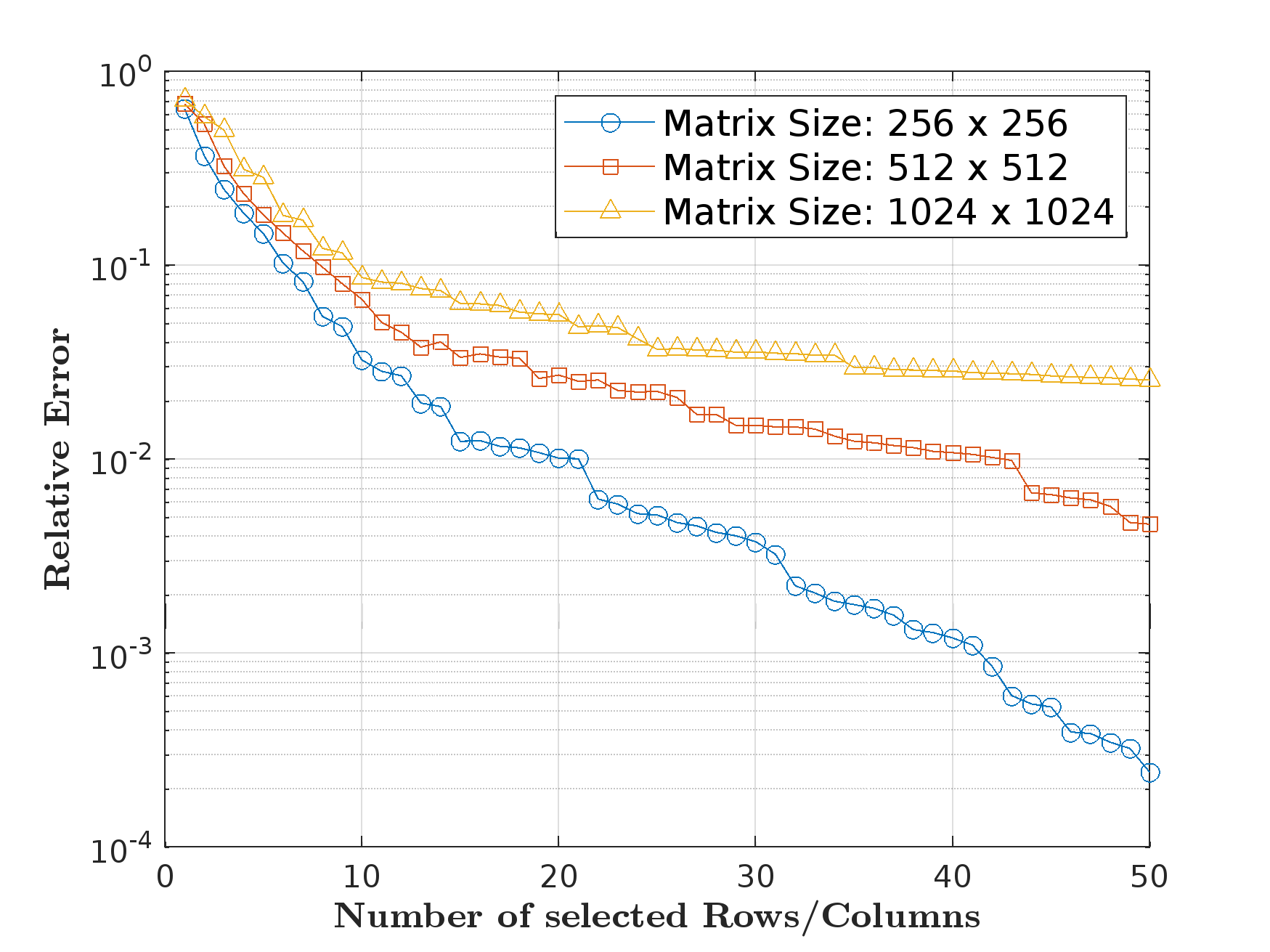}
\caption{Approximation error of \Cref{alg:Algorithm_CUR} for Hilbert matrices. The x-axis represents the number of selected rows/columns, and the y-axis shows the Frobenius norm relative error.}
\label{fig:hilbert_error}
\end{figure}

Hilbert matrices are well-known for their ill-conditioned nature and high sensitivity to numerical approximations~\cite{morhavc1995algorithm}. The experiments assess the approximation error of \Cref{alg:Algorithm_CUR} applied to Hilbert matrices of varying dimensions (\(256 \times 256\), \(512 \times 512\), and \(1024 \times 1024\)). For each matrix size, the number of selected rows and columns varied from \(1\) to \(20\). \Cref{fig:hilbert_error} shows the relationship between the number of selected rows and columns and the approximation error. 
The experimental results, illustrated in \Cref{fig:hilbert_error}, demonstrate the efficacy of \Cref{alg:Algorithm_CUR} in approximating matrices of varying sizes (\(256 \times 256\), \(512 \times 512\), and \(1024 \times 1024\)) by examining the relationship between the number of selected rows/columns and the Frobenius norm relative error. \Cref{fig:hilbert_error} reflects the algorithm's ability to effectively capture the dominant features of the matrix.

\paragraph{Experiments on Synthetic Low-Rank Matrices}

To further validate the algorithm, we conducted tests on synthetic low-rank matrices with controlled rank structures. These matrices were constructed by generating two matrices, \(A\) and \(B\), of sizes \(n \times r\) and \(r \times n\), respectively, with entries defined as:
\[
A(i, j) = \frac{1}{i + j}, \quad B(i, j)= \frac{1}{i + j},
\]
where \(i\) and \(j\) denote the row and column indices. The synthetic low-rank matrix \(H\) is generated by
\(
H = A \cdot B.
\)
The experiments are performed on the generated matrix with sizes \(n = 256, 512, 1024\) and rank value \(r = 5,10,15,20\).

\begin{figure}[htbp]
    \centering
    % First Row
    \begin{subfigure}[b]{0.48\textwidth}
        \includegraphics[width=\textwidth]{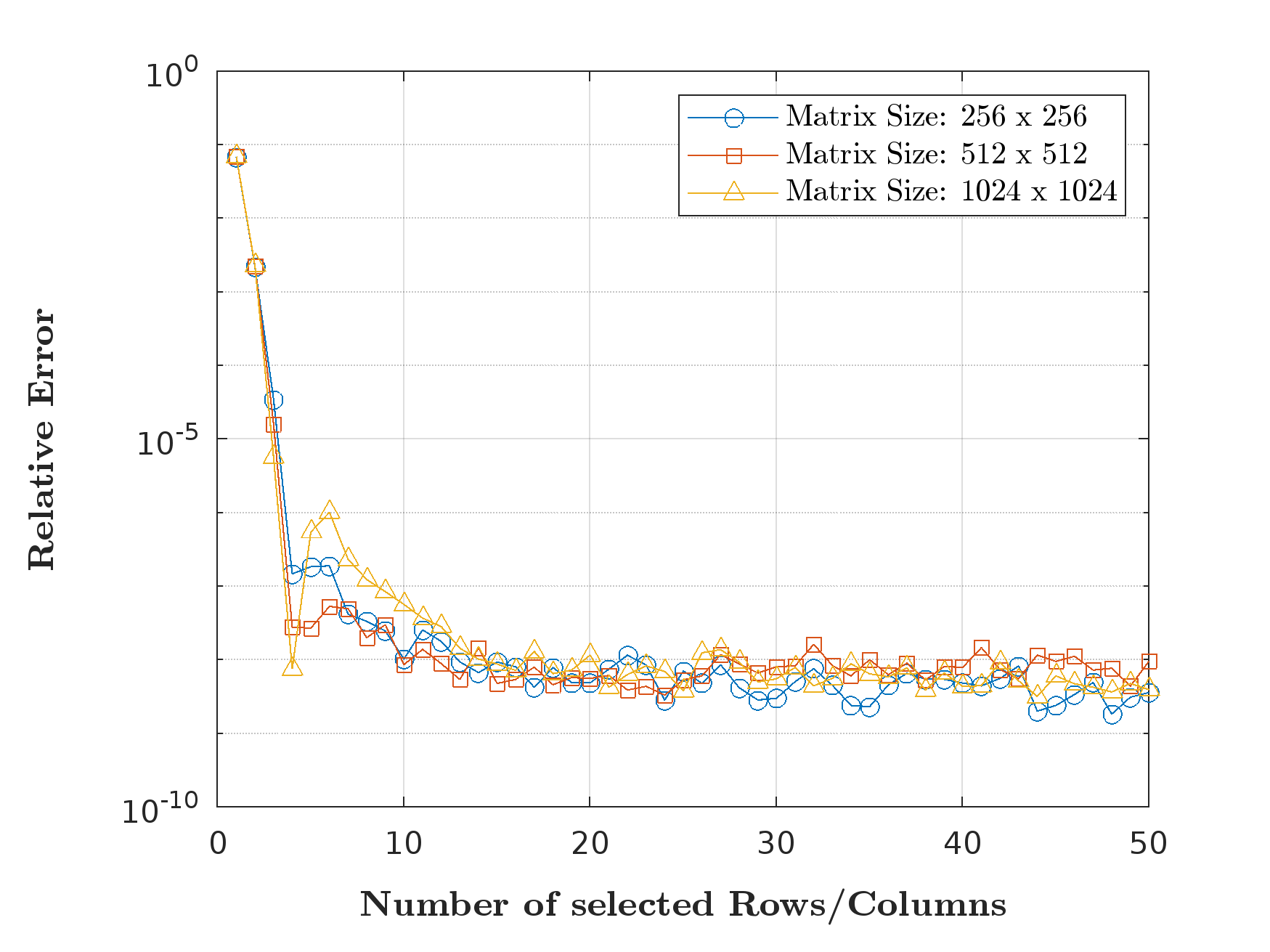}
        \caption{Rank = 5}
        \label{fig:elative_error_5}
    \end{subfigure}
    \hfill
    \begin{subfigure}[b]{0.48\textwidth}
        \includegraphics[width=\textwidth]{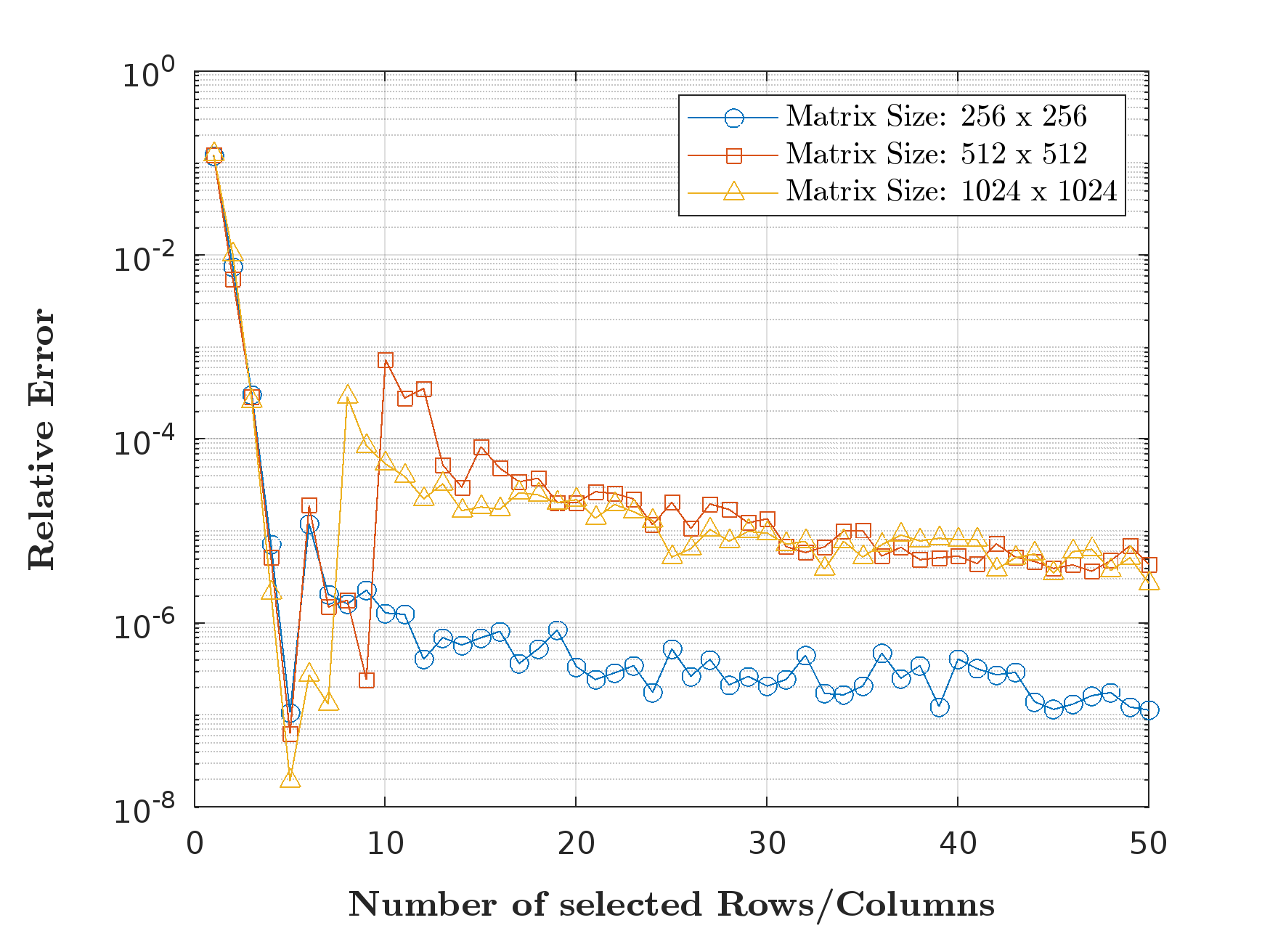}
        \caption{Rank = 10}
        \label{fig:relative_error_10}
    \end{subfigure}

    % Second Row
    \begin{subfigure}[b]{0.48\textwidth}
        \includegraphics[width=\textwidth]{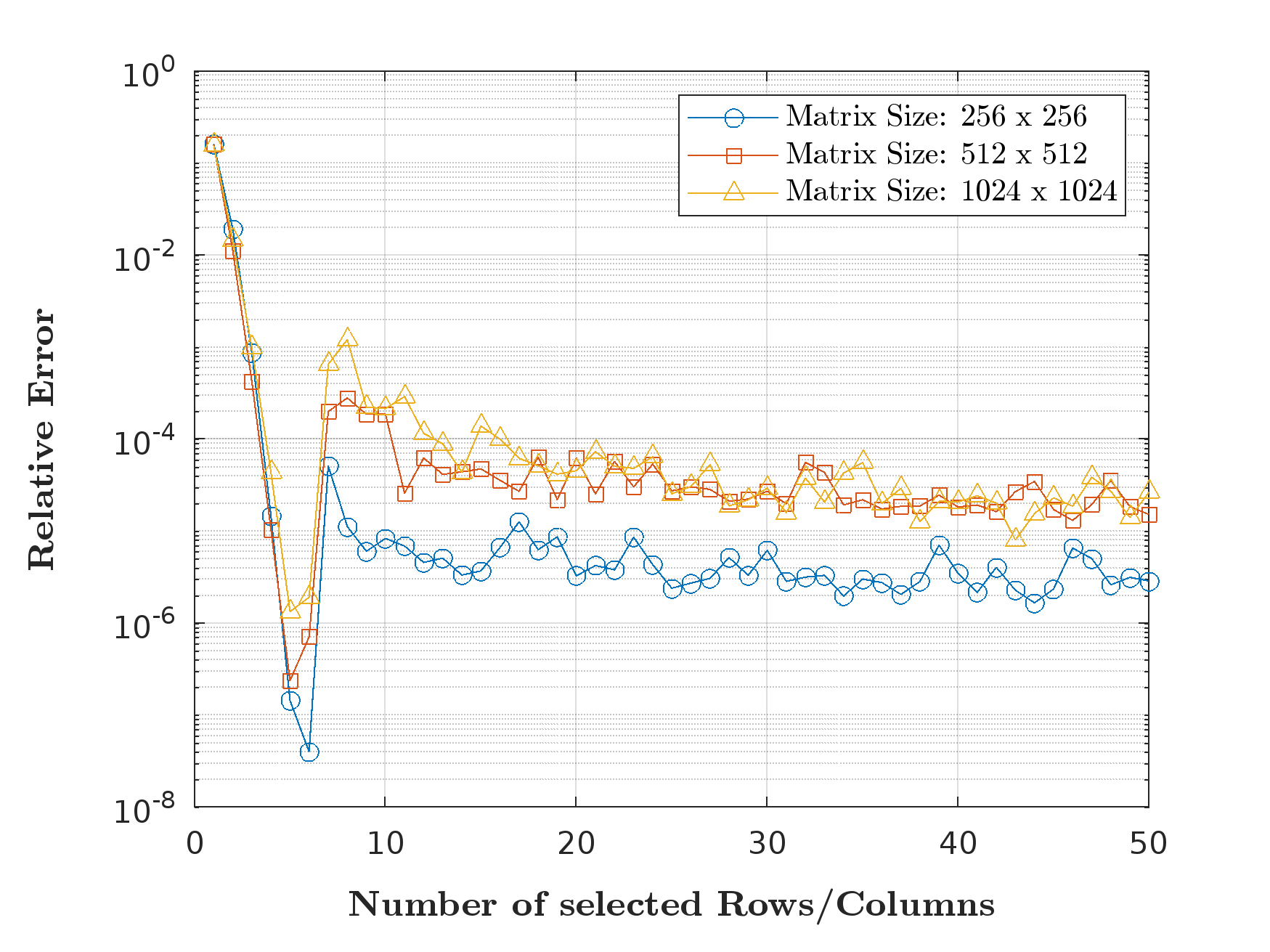}
        \caption{Rank = 15}
        \label{fig:relative_error_15}
    \end{subfigure}
    \hfill
    \begin{subfigure}[b]{0.48\textwidth}
        \includegraphics[width=\textwidth]{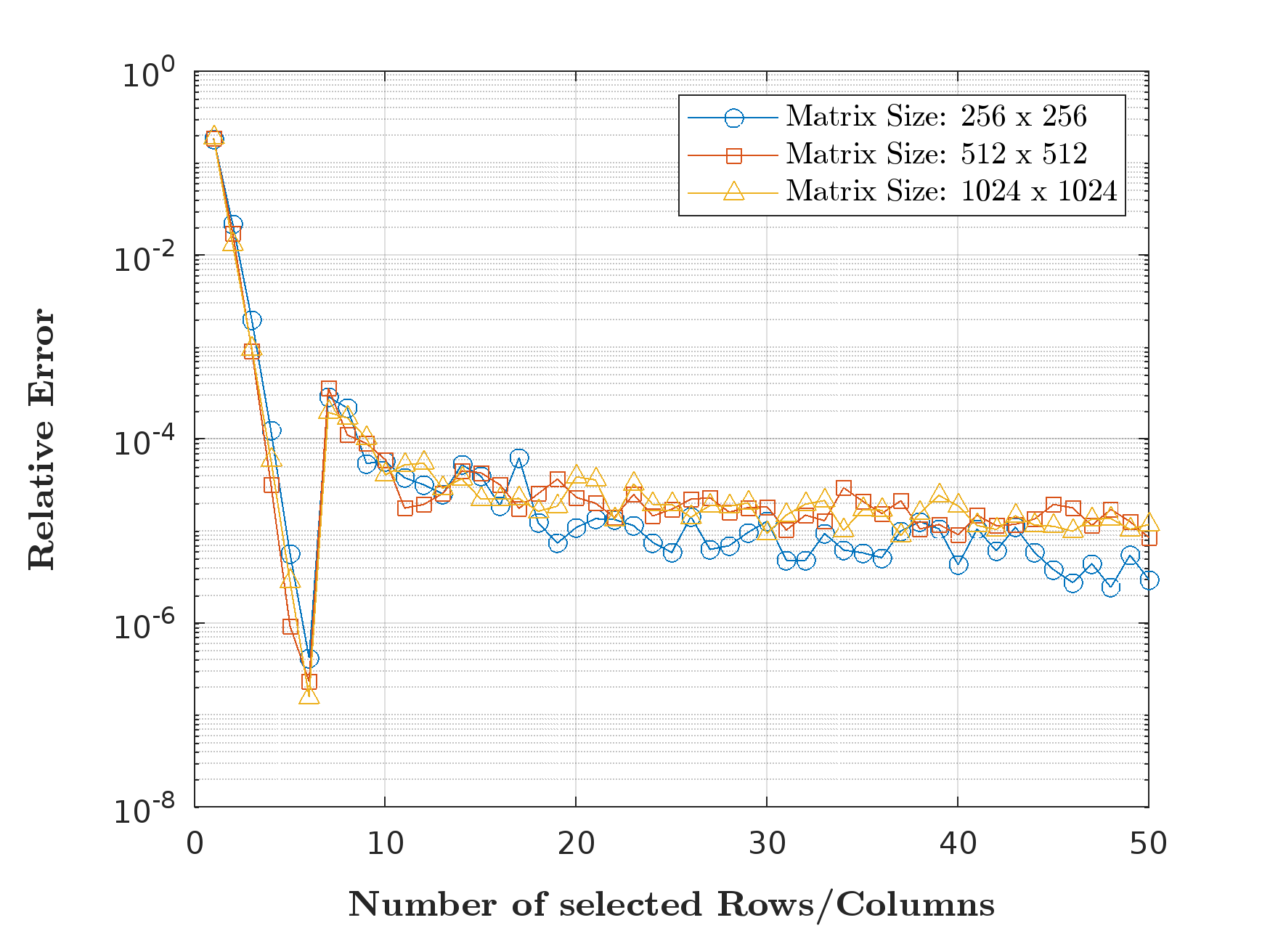}
        \caption{Rank = 20}
        \label{fig:relative_error_20}
    \end{subfigure}

    \caption{Relative Errors for Different Ranks}
    \label{fig:relative_error_all}
\end{figure}

The results reveal that the \Cref{alg:Algorithm_CUR} achieves near-perfect approximation for low-rank matrices even with a small number of selected rows and columns. The relative error decreases sharply as \(r\) increases and stabilizes at near-zero levels, reflecting the algorithm's effectiveness in reconstructing matrices with well-defined rank structures.

\subsection{Scalability Analysis of \Cref{alg:Algorithm_CUR}}

In this section, we conduct a scalability analysis of \Cref{alg:Algorithm_CUR}. The experiments focus on matrices of size \(16384 \times 16384\). We choose \Cref{alg:aca} as sequential benchmark algorithm 
\begin{figure}[H]
    \centering
    \includegraphics[width=1\linewidth]{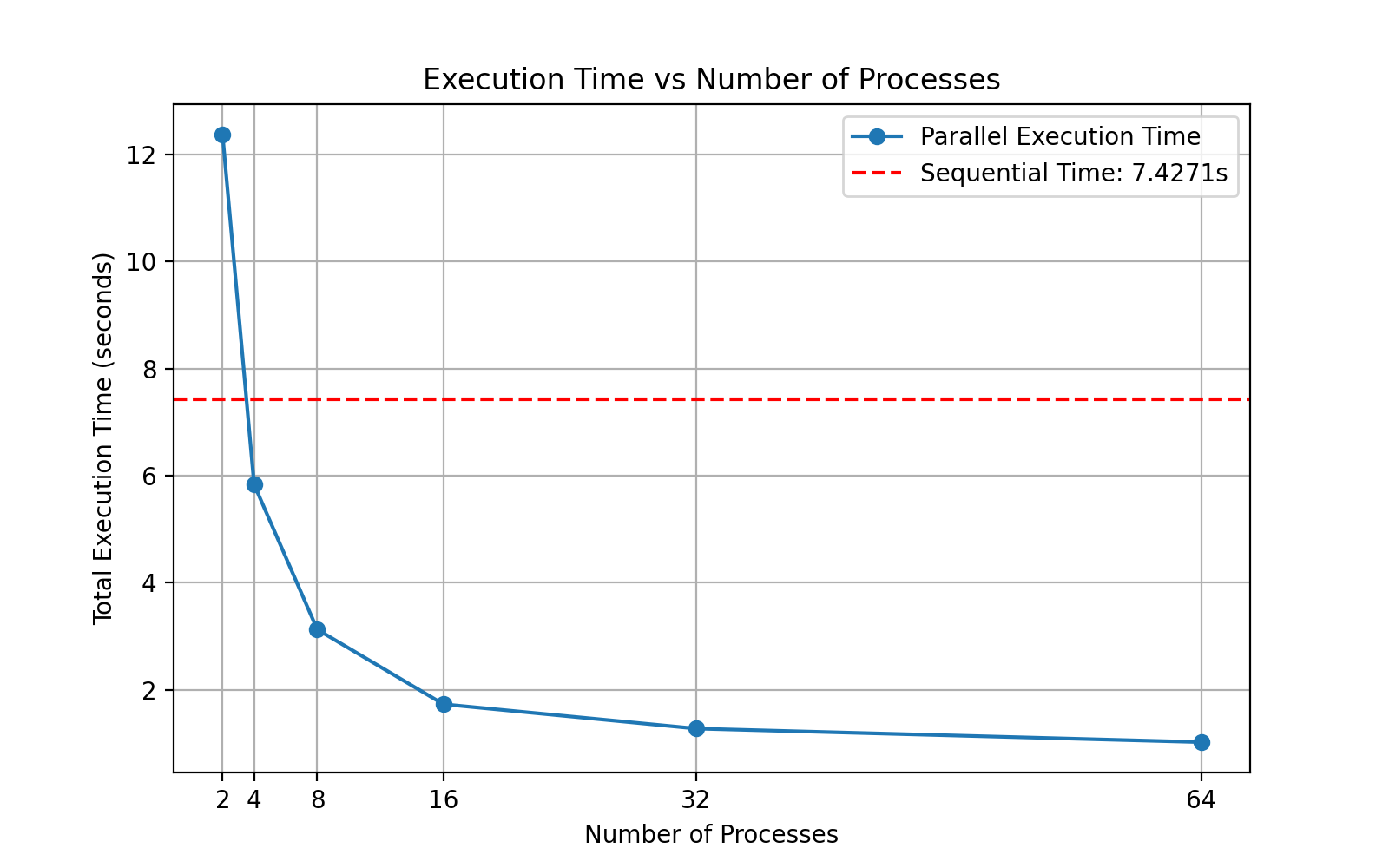}
    \caption{Caption}
    \label{fig:time}
\end{figure}
\Cref{fig:time} displays the relationship between the number of processes and the total execution time for the \Cref{alg:Algorithm_CUR}. As the number of processes increases from 2 to 64, the execution time exhibits a decreasing trend, demonstrating the expected performance improvement with parallelization. Notably, the execution time reduces from 12.3733 seconds for two process to 1.02225 seconds when employing 64 processes. However, the performance gain is not strictly linear; the curve gradually flattens as the number of processes increases beyond 16, suggesting diminishing returns due to parallel overhead and potential communication costs. A red dashed line is plotted to indicate the sequential baseline time of 7.4271 seconds, which serves as a reference point. The data highlights that parallel execution with more than 4 processes consistently outperforms the sequential baseline. These results suggest that parallelization is effective in accelerating MaxVol computations, but the efficiency gain saturates at higher process counts, likely due to inter-process communication and synchronization overhead.

\section{Conclusion}
In conclusion, we have presented a Scalable Binary CUR Low-Rank Approximation Algorithm. This approach is designed to achieve CUR decomposition of large-scale matrices by deterministiclly selecting rows and columns in parallel. By leveraging the parallel processing capabilities of modern multi-core architectures, our approach has shown modest improvements in computational efficiency. 

\bibliographystyle{abbrv}
\bibliography{reference}

\end{document}